\documentclass[psfig,a4paper,11pt,reqno]{amsart}
\usepackage{amsfonts,hyperref}
\usepackage{amssymb,amsmath,amsbsy,array,latexsym,color}

 \theoremstyle{plain}
\begingroup
\newtheorem{theorem}    {Theorem}                       [section]
\newtheorem{lemma}      [theorem]        {Lemma}
\newtheorem{proposition}[theorem]        {Proposition}
\newtheorem{corollary}  [theorem]        {Corollary}
\newtheorem{example}    [theorem]        {Example}
\theoremstyle{definition}
\newtheorem{remark}     [theorem]        {Remark}
\newtheorem{definition} [theorem]        {Definition}

\numberwithin{equation}                                 {section}

\title[sub-Riemannian objects in hypersurfaces]{On some sub-Riemannian objects in hypersurfaces
of sub-Riemannian manifolds}
\author[K. H. Tan, X. P. Yang]{Kang-Hai Tan, Xiao-Ping Yang}
\address{Department of Applied Mathematics, Nanjing
University of Science and Technology, 210094, Nanjing, Jiangsu, P.
R. China}

\email{tankanghai2000@yahoo.com.cn}\email{yangxp@mail.njust.edu.cn}

\begin{document}

\begin{abstract}
 \small{We study some sub-Riemannian objects (such as
horizontal connectivity, horizontal connection, horizontal tangent
plane, horizontal mean curvature) in hypersurfaces of sub-Riemannian
manifolds. We prove that if a connected hypersurface in a contact
manifold of dimension more than three is noncharacteristic or with
isolated characteristic points, then given two points, there exists
at least one piecewise smooth horizontal curve in this hypersurface
connecting them. In any sub-Riemannian manifold, we obtain the
sub-Riemannian version of the fundamental theorem of Riemannian
geometry: there exists a unique nonholonomic connection which is
completely determined by the sub-Riemannian structure and a
complement of the horizontal bundle, is ``symmetric'' and is
compatible with the sub-Riemannian metric. We use this nonholonomic
connection to study horizontal mean curvature of hypersurfaces.
 }
\end{abstract}
\maketitle 

\section{Introduction}
Recently there is an explosion of interest in the theory of
sub-Riemannian manifolds (or Carnot-Carath\'{e}odory spaces in
general), and in the ramification of this subject into analysis
and geometry, see e.g. \cite{Be}, \cite{DGN}, \cite{DGN1},
\cite{FSSC2}, \cite{FSSC3}, \cite{FSSC4}, \cite{FSSC5}, \cite{GN},
\cite{Gromov}, \cite{Magnani} and \cite{Pa}. We recall a
sub-Riemannian manifold $(M,\Delta,g_{c})$  is a smooth manifold
$M$ with a distribution $\Delta$ (a subbundle of the tangent
bundle $TM$) which is endowed with a fiberwise inner product
$g_{c}$ (usually called sub-Riemannian metric). $g_{c}$ is usually
realized as the restriction on $\Delta$ of some Riemannian metric
$g$ on $TM$. Carnot groups are particularly interesting
sub-Riemannian manifolds. Roughly speaking, the sub-Riemannian
geometry of $(M, \Delta,g_{c})$ is the geometry determined by the
sub-Riemannian structure $(\Delta,g_{c})$ which yields
Carnot-Carath\'{e}odory distance if $\Delta$ satisfies the Chow's
condition and $M$ is connected, and should be independent of the
choice of the Riemannian metric $g$ which is the extension of
$g_{c}$ to $TM$, although it often interacts with the Riemannian
geometry of $(M,g)$.

It is natural to study the geometry of hypersurfaces (or
submanifolds of codimension more than one) in sub-Riemannian
manifolds not only from the geometric viewpoint (see \cite{DGN},
\cite{GN}, \cite{GP} and \cite{Pa} for the program of
sub-Riemannian minimal surfaces in Carnot groups), but also from
the viewpoint of analysis (see \cite{FSSC2}-\cite{FSSC5},
\cite{Magnani} and \cite{Pa1} for the developing  geometric
measure theory in the setting of sub-Riemannian geometry). To
develop a theory of sub-Riemannian minimal surfaces, a
sub-Riemannian counterpart of the notion of the usual mean
curvature on submanifolds should be laid down. In \cite{Pa},
\cite{DGN} (see also \cite{GP}) an analysis definition of the
notion of horizontal mean curvature for noncharacteristic
hypersurfaces in Carnot groups has been given. In \cite{TY} we
have formulated a geometric definition of the horizontal mean
curvature by using the notion of horizontal connection. It turns
out that our definition coincides with the analysis definition in
the case of Carnot groups. The geometric definition is valid in
general sub-Riemannian manifolds.

   In this paper we continue our study on geometry and calculus of hypersurfaces
in sub-Riemannian manifolds. The notion of the horizontal tangent
planes of smooth noncharacteristic (see Definition
\ref{characteristic}) hypersurfaces plays an important role in the
story of \cite{TY}. The horizontal tangent plane
$T_{x}^{\mathcal{H}}S$ at a point $x$ in a smooth
noncharacteristic hypersurface $S$ is defined as the $k-1$
(assuming $\Delta$ is $k$ dimensional) dimensional linear space
 such that
$$
\Delta_{x}=n^{\mathcal{H}}(x)\bigoplus^{g_{c}}
T_{x}^{\mathcal{H}}S
$$
where $\Delta_{x}$ is the fiber of $\Delta$ through $x$,
$n^{\mathcal{H}}(x)$ the horizontal normal of $S$ at $x$ and
$\stackrel{g_{c}}{\bigoplus}$ denotes the orthogonal decomposition
with respect to $g_{c}$. Since $n^{\mathcal{H}}(x)$ is by
definition the projection onto $\Delta_{x}$ of the Riemannian
normal $n(x)$ computed with respect to $g$ (any orthogonal
extension of $g_{c}$), we have
$$
T_{x}^{\mathcal{H}}S=T_{x}S\bigcap \Delta_{x}.
$$
We call
$$
T^{\mathcal{H}}S:=\bigcup_{x\in{S}}T_{x}^{\mathcal{H}}S
$$
the horizontal tangent bundle. In general, if $S$ possesses
characteristic points, one also can define
$T^{\mathcal{H}}S=TS\bigcap \Delta$. Note that $T^{\mathcal{H}}S$
is independent of the orthogonal extensions of $g_{c}$. A very
interesting question is that under what condition
$T^{\mathcal{H}}S$ satisfies the Chow's condition with respect to
$TS$, that is, the subbundle $T^{\mathcal{H}}S$ (if it is a
subbundle, say) together with all its Lie commutators spans $TS$
(we assume $\Delta$ satisfies the Chow's condition). This question
is relevant to another two: whether there exist sufficiently many
horizontal curves in $S$ and whether one can define a
Carnot-Carath\'{e}odory metric with respect to $T^{\mathcal{H}}S$
with an induced metric $g^{\mathcal{H}}_{c}$ (from $g_{c}$)? We
point out that in his seminal paper \cite{Gromov} Gromov
discussed, among other things about geometry of hypersurfaces in
sub-Riemannian manifolds, the relationship between the restricted
Carnot-Carath\'{e}odory metric on $S$ and the
Carnot-Carath\'{e}odory metric induced by the structure
$(T^{\mathcal{H}}S, g^{\mathcal{H}}_{c})$. He claimed that the two
metrics are Lipschitz equivalent in some special cases, e.g.
contact manifolds of dimension more than three, see page 104 and
page 173 in \cite{Gromov}. But he did not seem to give conditions
to guarantee that the Carnot-Carath\'{e}odory metric induced by
the structure $(T^{\mathcal{H}}S, g^{\mathcal{H}}_{c})$ can be
well defined. One of main results of this paper is to give the
first result in this regard. We have
\begin{theorem}\label{Main1}
Let $M$ be a smooth orientable contact manifold of dimension
$2n+1$ with a contact form $\eta$ and a Riemannian metric $g$ and
let $\Delta:=\ker(\eta)$ be the canonical distribution. Let $S$ be
a smooth connected hypersurface in $M$. If $n>1$ and there does
not exist immersed submanifold contained in the set $\Sigma_{S}$
of all characteristic points in $S$, then for any points
$p,q\in{S}$, there exists a piecewise smooth curve
$\gamma(t),t\in[a,b]$ connecting $p,q$ such that
$\dot{\gamma}\in{T_{\gamma(t)}^{\mathcal{H}}S}$. In particular, if
$S$ is noncharacteristic, $T^{\mathcal{H}}S$ as a subbundle of
$TS$ satisfies the Chow's condition, that is, $T^{\mathcal{H}}S$
together with all its Lie commutators span $TS$.
\end{theorem}
Theorem \ref{Main1} follows from Sussmann's famous Orbit Theorem
and the contact structure of $M$. A trivial example shows that the
condition of dimension more than three is necessary.\vskip10pt

The next point we concentrate on is that whether there exists an
 intrinsic nonholonomic connection $D$ in $(M,\Delta,g_{c})$.
 A natural candidate is
the one $\widetilde{D}$ obtained by projecting to $\Delta$ the
Levi-Civita connection $\nabla$ with respect to some orthogonal
extension $g$ of $g_{c}$. We prove that $\widetilde{D}$ is
independent of the choice of orthogonal extension of $g_{c}$.
Moreover we have
\begin{theorem}\label{Main2} Let $\{X_{1},\cdots,X_{k}\}$ be an
orthonormal basis of $\Delta$. We define $D$ by
\begin{equation}\label{definition}
D_{U}V=\sum_{i=1}^{k}U(V^{i})X_{i}+\sum_{i=1}^{k}\sum_{j=1}^{k}\sum_{l=1}^{k}U^{j}V^{i}\Gamma_{ij}^{l}X_{l}
\textrm{ for any }U=\sum_{j=1}^{k}U^{j}X_{j},\,
V=\sum_{i=1}^{k}V^{i}X_{i}\in{\Gamma(\Delta)}
\end{equation}
where
$\Gamma_{ij}^{l}=-\frac{1}{2}\{g_{c}(X_{i},[X_{j},X_{l}]^{\mathcal{H}})+
g_{c}(X_{j},[X_{l},X_{i}]^{\mathcal{H}})
-g_{c}(X_{l},[X_{i},X_{j}]^{\mathcal{H}})\}$ and $X^{\mathcal{H}}$
is understood as the projection of $X$ to $\Delta$ with respect to
the direct summation decomposition
$TM=\Delta\bigoplus\widetilde{\Delta}$ where $\widetilde{\Delta}$
is the distribution complementary to $\Delta$. Then $D$ is
independent of the choice of orthonormal basis of $\Delta$ and it
is the \textbf{unique} nonholonomic connection  satisfying
\begin{enumerate}
    \item\label{1'} $D_{U}V$ is $R-$linear in both arguments,
    \item\label{2'} $D_{U}V$ is $C^{\infty}(M)-$linear in the argument of
    $U$,
    \item\label{3'} the Leibniz rule holds:
    $$
    D_{U}(fV)=(Uf)V+fD_{U}V\textrm{ for any }
    f\in{C^{\infty}(M)}, U,V\in{\Gamma(\Delta)},
    $$
    \item $D$ is compatible with respect to $g_{c}$ , that
    is,
    \begin{equation}\label{4'}
   Ug_{c}(V,W)=g_{c}(D_{U}V,W)+g_{c}(U,D_{V}W)\textrm{ for any
    }U,V,W\in{\Gamma(\Delta)},
   \end{equation}
   \item the following ``symmetry'' property holds:
     \begin{equation}\label{symmetry2}
 D_{U}V-D_{V}U=[U,V]^{\mathcal{H}}\textrm{ for any
 }U,V\in{\Gamma(\Delta)}.
\end{equation}
\end{enumerate}
In particular we have $D=\widetilde{D}$.
\end{theorem}
Theorem \ref{Main2} is the counterpart of the existence and
uniqueness of the Levi-Civita connection in Riemannian geometry. But
the connection certainly depends on the splitting of $TM$.

We will use $D$ to study the horizontal mean curvature of
hypersurfaces in $(M,\Delta,g_{c})$. Let $D^{\top}$ be the tangent
horizontal connection on the horizontal tangent bundle
$T^{\mathcal{H}}S$ of a smooth noncharacteristic hypersurface $S$.
We define the horizontal mean curvature of $S$ as the trace of the
horizontal fundamental second form which is by definition a
bilinear map $\textrm{II}$ from $\Gamma(T^{\mathcal{H}}S)\times
\Gamma(T^{\mathcal{H}}S)$ to $N$:
$$
\textrm{II}(X,Y)=D_{X}Y-D^{\top}_{X}Y
$$
for any $X,Y\in{\Gamma(T^{\mathcal{H}}S)}$ where
$\Gamma(T^{\mathcal{H}}S)$ denotes the set of all smooth sections
of $T^{\mathcal{H}}S$ and $N$ is the horizontal normal. The
symmetry of $\textrm{II}$ follows from the symmetry property
\eqref{symmetry2} of $D$ and the definition of the horizontal
tangent plane. Since both $D$ and $T^{\mathcal{H}}S$ are
intrinsic, so is the horizontal mean curvature.

The paper is organized as follows. In Section \ref{sec:basic} we
collect some facts about sub-Riemannian manifolds which will be
used later, mainly to fix some notations. In Section
\ref{sec:objects} Theorem \ref{Main1} is proven after introducing
the notion of the horizontal tangent bundle. Roughly speaking, if
we project Riemannian objects onto the horizontal bundle, such as
Riemannian connection, normal vector and tangent bundle, then we
get corresponding sub-Riemannian analogues: horizontal connection,
horizontal normal vector and horizontal bundle. Theorem
\ref{Main2} is proven in Section \ref{sec:proof2}.

\section{Basic material on sub-Riemannian manifolds}\label{sec:basic}
Let $M$ be a smooth ($C^{\infty}$) manifold of dimension $m$ endowed
with a smooth distribution (called \textit{horizontal bundle})
$\Delta$ of dimension $k$ with $k<m$.\footnote{Note that this
imposes topological constraints on $M$, see \cite{Ko}} If we a prior
equip $\Delta$ with an inner product $g_{c}$ (called
\textit{sub-Riemannian metric}), we call $(M,\Delta,g_{c})$ a
\textit{sub-Riemannian manifold} with the \textit{sub-Riemannian
structure} $(\Delta,g_{c})$. Let $\{X_{1},\cdots,X_{k}\}$ be an
orthonormal local basis of $\Delta$. A piecewise smooth curve
$\gamma(t),t\in{[a,b]}$ in $M$ is horizontal if
$\dot{\gamma(t)}\in{\Delta_{\gamma(t)}}$ a.e. $t\in{[a,b]}$. The
\textit{length} $\ell(\gamma)$ of the horizontal curve
$\gamma(t),t\in{[a,b]}$ is the integral
$\int_{a}^{b}g_{c}(\dot{\gamma}(t),\dot{\gamma}(t))dt$. Denote by
$L_{i}$ the set of all vector fields spanned by all commutators of
$X_{j}$'s of order $\leq i$ and let $L_{i}(p)$ be the subspace of
evaluations  at $p$ of all vector fields in $L_{i}$. We call
$\Delta$ satisfies the \textit{Chow's condition}\footnote{In the
subelliptic theory, Chow's condition is also called H\"{o}rmander's
condition.} if for any $p\in{M}$, there exists an integer $r(p)$
such that $L_{r(p)}(p)=T_{p}M$ (the tangent space of $M$ at $p$). If
$M$ is connected and $\Delta$ satisfies the Chow's condition, the
Chow connectivity theorem asserts that there exists at least one
piecewise smooth horizontal curve connecting two given points (see
\cite{Chow}, \cite{Be} or \cite{Gromov}), and thus $(\Delta,g_{c})$
yields a metric (called \textit{Carnot-Carath\'eodory metric})
$d_{c}$ by letting $d_{c}(p,q)$ as the infimum among the lengths of
all horizontal curves joining $p$ to $q$.

$\Delta$ is \textit{equiregular} if the dimension of $L_{i}(p)$
does not depends on $p$ for any $i$, that is, the tangent bundle
$TM$ is filtered by smooth subbundles
$$
\Delta=L_{1}\subset L_{2}\subset\cdots\subset L_{r}=TM
$$
where $r$ is called the degree of $\Delta$. Note that we always
can extend $g_{c}$ to a Riemannian metric $g$ in $M$ such that
$TM$ can be $g$-orthogonally decomposed as
$TM=\Delta\stackrel{g}{\bigoplus}\widetilde{\Delta}$ where
$\widetilde{\Delta}$ is the distribution complementary to
$\Delta$. We call such $g$ an \textit{orthogonal extension} of
$g_{c}$. Obviously the orthogonal extension of $g_{c}$ is not
unique in general. We will use $\Gamma(\Delta)$ to denote the set
of all smooth sections of $\Delta$.
\begin{example}[\textbf{Carnot groups}] The most interesting models of
sub-Riemannian manifolds are Carnot groups (called also stratified
groups). A \textit{Carnot group} $G$ is a connected, simply
connected Lie group whose Lie algebra $\mathcal{G}$ admits the
grading $\mathcal{G}=V_{1}\bigoplus\cdots\bigoplus V_{l}$, with
$[V_{1},V_{i}]=V_{i+1}$, for any $1\leq i\leq l-1$ and
$[V_{1},V_{l}]={0}$ (the integer $l$ is called the step of $G$).
Let $\{e_{1},\cdots,e_{n}\}$ be a basis of $\mathcal{G}$ with
$n=\sum_{i=1}^{l}\dim(V_{i})$. Let $X_{i}(g)=(L_{g})_{*}e_{i}$ for
$i=1,\cdots,k:=\dim(V_{1})$ where $(L_{g})_{*}$ is the
differential of the left translation
$L_{g}(g^{\prime})=gg^{\prime}$ and let
$Y_{i}(g)=(L_{g})_{*}e_{i+k}$ for $i=1,\cdots,n-k$. We call the
system of left-invariant vector fields
$\Delta:=V_{1}=\textrm{span}\{X_{1},\cdots,X_{k}\}$ the
\textit{horizontal bundle} of $G$. If we equip $\Delta$ an inner
product $g_{c}$ such that $\{X_{1},\cdots,X_{k}\}$ is an
orthonormal basis of $\Delta$, $(G,\Delta,g_{c})$ is an
equiregular sub-Riemannian manifold. In $(G,\Delta,g_{c})$,
$d_{c}$ is invariant with respect to left translation, that is
$d_{c}(p_{0}p,p_{0}q)=d_{c}(p,q)$ for any $p_{0},p,q\in{G}$, and
is 1-homogeneous with respect to the natural dilations, that is
$d_{c}(\delta_{s}p,\delta_{s}q)=sd_{c}(p,q)$ for any
$s>0,p,q\in{G}$, where
$\delta_{s}p=\exp(\sum_{i=1}^{l}s^{i}\xi_{i})$ for
$p=\exp(\sum_{i=1}^{l}\xi_{i}), \xi_{i}\in{V_{i}}$.
\end{example}
\begin{example}[\textbf{contact manifolds}]Contact manifolds are also
interesting sub-Riemannian manifolds not necessarily with group
structure. A contact manifold is a smooth (connected) manifold $M$
with a contact structure which is by definition a codimension one
distribution $\Delta\subset TM$ with non-degenerate curvature form
$\omega:\Delta\wedge\Delta\rightarrow TM/\Delta$ which is defined
as follows: first represent $\Delta$ locally as the kernel of a
1-form, say $\eta$ on $M$, then identify $TM/\Delta$ with the
trivial line bundle and finally define $\omega$ as $d\eta|\Delta$.
$\eta$ is called a contact form of $M$ (not unique). If $M$ is
orientable, $\eta$ and $\omega$ can be globally defined. The
non-degeneracy of $\omega$ makes the dimension of $\Delta$ be even
and so the dimension of $M$ is odd, say $2n+1$. Note that the
non-vanishing of $\omega$ on $\Delta$ makes the commutators of
degree $\leq2$ span $TM$ and thus $(M,\Delta,g_{c})$ where $g_{c}$
is the restriction on $\Delta$ of some Riemannian metric $g$ on
$M$ is an equiregular sub-Riemannian manifold.
\end{example}
\begin{example}[\textbf{Heisenberg group}]\label{Heisenberg} Heisenberg group $\mathbb{H}^{n}$ as a
representative both in the class of Carnot groups and in the class
of contact manifolds is of paramount importance, and is worthy of
being paid more attention to. The underlying manifold of this Lie
group is simply $\mathbb{R}^{2n+1}$, with the noncommutative group
law
$$
p\cdot
p^{\prime}=(x,y,t).(x^{\prime},y^\prime,t^{\prime})=(x+x^{\prime},y+y^{\prime},t+t^{\prime}+
2(<x^{\prime},y>-<x,y^{\prime}>))
$$
 where we have let $x,
x^{\prime}, y, y^{\prime}\in{\mathbb{R}^{n}}, t,
t^{\prime}\in{\mathbb{R}}$. A simple computation shows that
$$\begin{array}{cccl}
  X_{j}(p)&\stackrel{\textrm{def}}{=}&(L_{p})_{*}(\frac{\partial}{\partial
x_{j}})=&\frac{\partial}{\partial
x_{j}}+2y_{j}\frac{\partial}{\partial t},\quad j=1,\cdots,n,\\
  X_{n+j}(p)&\overset{\textrm{def}}{=} &(L_{p})_{*}(\frac{\partial}{\partial
y_{j}})=&\frac{\partial}{\partial
y_{j}}-2x_{j}\frac{\partial}{\partial t},\quad j=1,\cdots,n,\\
T(p)&\overset{\textrm{def}}{=}&(L_{p})_{*}(\frac{\partial}{\partial
t})=&\frac{\partial}{\partial t}\quad
\end{array}$$ for any $p=(x,y,t)$ in $\mathbb{H}^{n}$.
We note that
$$
\begin{array}{c}
[X_{j},X_{n+k}]=-4T\delta_{jk},\quad j,k=1,\cdots,n, \\
\textrm{and all other commutators are trivial},\end{array}
$$
therefore the vector fields $X=\{X_{1},\cdots,X_{2n}\}$ constitute
a basis of the Lie algebra
$\flat_{n}=\mathbb{R}^{2n+1}=V_{1}\oplus V_{2}$, where
$V_{1}=\mathbb{R}^{2n}\times\{0\}_{t},V_{2}=\{0\}_{x,y}\times\mathbb{R}.$
Note that the horizontal bundle
$\Delta=\textrm{span}\{X_{1},\cdots,X_{2n}\}$ is the kernel of the
1-form
$\eta=\frac{1}{4}dt+\frac{1}{2}\sum_{i=1}^{n}(x_{i}dy_{i}-y_{i}dx_{i})$
and the curvature form $\omega=d\eta=\sum_{i=1}^{n}dx_{i}\wedge
dy_{i}$ is the standard symplectic form in $\mathbb{R}^{2n}$. Thus
a smooth curve $\gamma(s)=(x(s),y(s),t(s)):[a,b]\rightarrow
\mathbb{H}^{n}$ is horizontal if and only if
\begin{equation}\label{horizontal}
2\dot{t}(s)=\sum_{i=1}^{n}y_{i}(s)\dot{x_{i}}(s)-x_{i}(s)\dot{y}_{i}(s)\quad\textrm{
for any }s\in[a,b].
\end{equation}
\end{example}
For the theory of sub-Riemannian geodesics we refer to the book
\cite{Montgomery} and references therein, and see in particular
\cite{Gromov} for a comprehensive treatment of the geometry (more
than sub-Riemannian geodesics) in sub-Riemannian manifolds, also for
many potential research directions.

\section{The horizontal
tangent bundle and horizontal connectivity in
hypersurfaces}\label{sec:objects}

 Let $(M,\Delta,g_{c})$ be a sub-Riemannian manifold. \textbf{In this section we always assume
 $M$ is connected and $\Delta$ satisfies the Chow's condition.} By $S$ we always mean a
smooth hypersurface (i.e. an embedded submanifold of codimension
1) in $M$. Let $g$ be any orthogonal extension of $g_{c}$.

\begin{definition}[\textbf{characteristic points}]\label{characteristic}A point $p\in{S}$ is a \textit{characteristic} point if
$\Delta_{p}\subset T_{p}S$. Let $\Sigma_{S}$ denote the set of all
characteristic points in $S$. If $\Sigma_{S}=\emptyset$, we call
$S$ is \textit{noncharacteristic}.
\end{definition}
Typically $S$ possesses characteristic points, see e.g.
\cite{DGN1}. But we have
\begin{proposition}\label{smallness}
     $\Sigma_{S}$ is a closed subset of $S$, and
    $\mathcal{H}^{m-1}(\Sigma_{S})=0$. Here $\mathcal{H}^{m-1}$
    denotes the $m-1$ dimensional Hausdorff measure with respect
    to the Riemannian metric $g$ (recall $M$ is $m$-dimensional).
\end{proposition}
The closedness follows from the smoothness of $S$ and the
smallness of $\Sigma_{S}$ is due to Derridj, see \cite{De1},
\cite{De2}. In the case of Carnot groups of step two, when the
smoothness of $S$ is weaker than $C^{\infty}$, say $C^{1,1}$, the
smallness of $\Sigma_{S}$ is obtained by Magnani, see
\cite{Magnani} or \cite{Magnani1}.

By Proposition \ref{smallness}, if $p\in{S\backslash\Sigma_{S}}$,
then there exists a neighborhood $\mathcal{U}$ of $p$, such that
$\mathcal{U}\cap S\subset S\backslash\Sigma_{S}$.
\begin{definition}[\textbf{horizontal normal}]\label{horizontal normal}
Let $n^{g}$ denote the Riemannian normal of $S$ with respect to
$g$. We define the \textit{horizontal normal}  $n^{\mathcal{H}}$
of $S$ as the projection of the Riemannian normal onto the
horizontal bundle, that is,
$$
n^{\mathcal{H}}=\sum_{i=1}^{k}g(n^{g},X_{i})X_{i}
$$
where $\{X_{1},\cdots,X_{k}\}$ is an orthonormal basis of
$\Delta$.
\end{definition}
It is easily seen that $p\in{S}$ is a characteristic point if and
only if $n^{\mathcal{H}}(p)=0$.
\begin{definition}[\textbf{horizontal tangent plane}]For $p\in{S}$,
we define
$$
T_{p}^{\mathcal{H}}S:=\{v\in{\Delta_{p}}:g_{c}(v,n^{\mathcal{H}}(p))=0\}
$$
as the \textit{horizontal tangent plane} of $S$ at $p$.
\end{definition}
From the definition we see that if $p\in{\Sigma_{S}}$ then
$T_{p}^{\mathcal{H}}S=\Delta_{p}$, and otherwise
$T_{p}^{\mathcal{H}}S$ is a $k-1$ dimensional subspace of
$\Delta_{p}$. Thus
$$
T^{\mathcal{H}}S:=\bigcup_{p\in{S}}T_{p}^{\mathcal{H}}S
$$
is a distribution of dimension $k-1$ if and only if $S$ is
noncharacteristic. The following proposition shows that
$T^{\mathcal{H}}S$ depends only on $\Delta$ and $S$.
\begin{proposition}The horizontal tangent bundle $T^{\mathcal{H}}S$
is intrinsic:
$$
T^{\mathcal{H}}S=\Delta\bigcap TS
$$
in the sense that $T_{p}^{\mathcal{H}}S=\Delta_{p}\bigcap T_{p}S$
for any $p\in{S}$.
\end{proposition}\label{horizontalbundle}
\proof Let $g$ be any orthogonal extension of $g_{c}$. For
$p\in{S}$, first let $v\in{T_{p}^{\mathcal{H}}S}$. Thus
$v\in{\Delta_{p}}$ and from Definition \ref{horizontal normal} we
have
$$
\begin{aligned}
g(v,n^{g}(p))&=g_{c}(v,n^{\mathcal{H}}(p))+g(v,n^{\prime}(p))\\
&=0+0=0,
\end{aligned}
$$
where $n^{\prime}(p)=n^{g}(p)-n^{\mathcal{H}}(p)$ is orthogonal to
$v$ since $g$ is an orthogonal extension. The last formulas
implies that $v\in{T_{p}S}$. So
$T_{p}^{\mathcal{H}}S\subset\Delta_{p}\bigcap T_{p}S$.
$\Delta_{p}\bigcap T_{p}S\subset T_{p}^{\mathcal{H}}S$ follows
from the same argument.
\endproof
\begin{remark}\label{normalintrinsic}
Therefore $T^{\mathcal{H}}S$ is the projection of $TS$ onto
$\Delta$. In the case of Carnot groups of step two, if $p\in{S}$
is not a characteristic point, $T_{p}^{\mathcal{H}}S$ has obvious
geometric meaning: $T^{\mathcal{H}}_{p}S$ is the projection on the
horizontal bundle of the Lie algebra of the tangent group which is
the blowup set of $S$ with respect to the natural dilations
\footnote{That is the limit set of
$S_{p,s}:=\delta_{\frac{1}{s}}(p^{-1}S)$ under suitable topology
as $s\rightarrow0$.}, see \cite{FSSC4}. This is the reason why we
call $T_{p}^{\mathcal{H}}S$ horizontal tangent plane.

If $p\in{S\backslash\Sigma_{S}}$, let
$\mathcal{V}(p)=\frac{n^{\mathcal{H}}(p)}{|n^{\mathcal{H}}(p)|}$
be the unit horizontal normal. Then by Proposition
\ref{horizontalbundle}, $\mathcal{V}$ is intrinsic: independent of
the choice of orthogonal extension, since
$\Delta_{p}=T_{p}^{\mathcal{H}}S\stackrel{g_{c}}{\bigoplus}
\mathcal{V}(p)$.
\end{remark}

 A natural question is that whether there exist sufficiently
many horizontal curves in hypersurfaces such that the intrinsic
Carnot-Carath\'{e}odory metric can be defined. We will not pursue
the general case. But we will prove that for contact manifolds of
dimension more than three, there exists at least one smooth
horizontal curve connecting two given points in a connected
hypersurface $S$ if $\Sigma_{S}$ does not contain any immersed
submanifold, and so the intrinsic Carnot-Carath\'{e}odory metric
can be defined. This is the case if $S$ is noncharacteristic or
with isolated characteristic points. Note that even in the contact
case, the horizontal connectivity in hypersurfaces is overlooked.
Some authors asserted that there are few horizontal curves in
hypersurfaces, even though these hypersurfaces are
noncharacteristic, see e.g. \cite{FSSC2} (p.485) where the authors
wrote: ``$\cdots$Notice however that a H-regular hypersurface
contains very few H-rectifiable curves; in particular we cannot
define a geodesic distance on a H-rectifiable
hypersurface$\cdots$''\footnote{Smooth noncharacteristic
hypersurfaces are, from the definition, H-regular hypersurfaces;
the inverse is not true in general. H-rectifiable curves are just
absolutely continuous horizontal curves.}. \vskip10pt

 To prove Theorem \ref{Main1}, we first
introduce the notion of the orbit of a family of vector fields.
\begin{definition}[\textbf{orbits of a family of vector fields}]Let $M$ be
a connected smooth manifold and let $\mathcal{F}$ be any family of
smooth vector fields globally defined on $M$. We define the orbit
of a point $p\in{M}$ of this family as the set of points of $M$
reachable piececwise by trajectories of vector fields in the
family, that is,
$$
\mathcal{O}_{p}:=\{\exp(t_{n}f_{k})\circ\cdots\circ\exp(t_{1}f_{1})\circ
p\mid
t_{i}\in{\mathbb{R}},f_{i}\in{\mathcal{F}},n\in{\mathbb{N}}\}
$$
where $\exp(tf)(p)$ denotes the flow of the vector field $f$
through $p$, i.e. the curve $\gamma(t)$ in $M$ such that
$$
\left\{
\begin{aligned}
&\dot{\gamma}(t)=f(\gamma(t))\\
&\gamma(0)=p
\end{aligned}
\right..
$$
Of course, if some of our vector fields are not complete then we
consider only such $t_{1},\cdots,t_{n}$ for which the above
expression has sense. It is clear that the relation:``$q$ belongs
to the orbit of $p$'' is an equivalence relation on $M$ and thus
$M$ is the disjoint union of orbits (equivalence classes).
\end{definition}
The following orbit theorem is due to Sussmann (also Nagano), see
\cite{Sussmann}.
\begin{theorem}[\textbf{Orbit Theorem, Nagano-Sussmann}]\label{orbit}Let
$\mathcal{F}$ be as above and let $p\in{M}$.
Then:\begin{enumerate}
    \item $\mathcal{O}_{p}$ is a connected immersion submanifold
    of $M$.
    \item $T_{q}\mathcal{O}_{p}=\textrm{span}\{(P_{*}f)(q)\mid
    P\in{\mathcal{P},f\in{\mathcal{F}}}\}, q\in{\mathcal{O}_{p}}$
    where $\mathcal{P}$ is the group of diffeomorphisms
    of $M$ generated by flows in $\mathcal{F}$:
    $$
    \mathcal{P}=\{\exp(t_{n}f_{k})\circ\cdots\circ\exp(t_{1}f_{1})
    \mid
    t_{i}\in{\mathbb{R}},f_{i}\in{\mathcal{F}},n\in{\mathbb{N}}\}\subset
    \textrm{Diff}(M)
    $$
and $P_{*}$ is the differential map of $P$.
    \end{enumerate}
\end{theorem}
\begin{remark}From Theorem \ref{orbit}, two simple but
very useful observations are in order.
\begin{enumerate}
    \item First of all, if $f\in{\mathcal{F}}$, then
    $f(q)\in{T_{q}\mathcal{O}_{p}}$ for all
    $q\in{\mathcal{O}_{p}}$. Indeed, the trajectory $\exp(tf)(q)$
    belongs to the orbit $\mathcal{O}_{p}$, thus its velocity
    vector $f(q)$ is in the tangent space $T_{q}\mathcal{O}_{p}$.
    \item Further, if $f_{1},f_{2}\in{\mathcal{F}}$, then
    $[f_{1},f_{2}](q)\in{T_{q}\mathcal{O}_{p}}$ for all
    $q\in{\mathcal{O}_{p}}$. This follows since the vector
    $[f_{1},f_{2}](q)$ is tangent to the trajectory
    $\exp(-tf_{2})\circ\exp(-tf_{1})\circ\exp(tf_{2})\circ\exp(tf_{1})(q)\in{\mathcal{O}_{p}}$.
    We go on and consider Lie brackets of arbitrarily high order
    $$
    [f_{1},[\cdots[f_{n-1},f_{n}]\cdots]](q)
    $$
    as tangent vectors to $\mathcal{O}_{p}$ if
    $f_{i}\in{\mathcal{F}}$ and $q\in{\mathcal{O}_{p}}$.
    \end{enumerate}
These considerations can be summarized in terms of Lie algebra of
vector fields generated by $\mathcal{F}$:
$$
\textrm{Lie}\mathcal{F}:=\textrm{span}\{[f_{1},[\cdots[f_{n-1},f_{n}]\cdots]]\mid
f_{i}\in{\mathcal{F}},n\in{\mathbb{N}}\}\subset (\Gamma(TM)),
$$
and its evaluation at a point $q\in{M}$:
$$
\textrm{Lie}_{q}\mathcal{F}=\{V(q)\mid
V\in{\textrm{Lie}\mathcal{F}}\}\subset T_{q}M.
$$
\end{remark}
We obtain the following statement.
\begin{corollary}\label{contained}
$$
\textrm{Lie}_{q}\mathcal{F}\subset T_{q}\mathcal{O}_{p}
$$
for all $q\in{\mathcal{O}_{p}}$.
\end{corollary}
We note that the Chow connectivity theorem follows immediately
from Corollary \ref{contained}.

Now we return to consider the horizontal connectivity in
hypersurfaces in sub-Riemannian manifolds. We will need the notion
of a horizontal immersed submanifold with respect to a
distribution.
\begin{definition}
Let $M$ be a smooth manifold and let $\Delta$ be a distribution on
$M$. An immersed submanifold $i:N\rightarrow M$ of $M$ is
horizontal with respect to $\Delta$ if $i_{*}(T_{p}N)\subset
\Delta_{p}$ for any $p\in{N}$.
\end{definition}

For codimension one distribution we have
\begin{theorem}\label{codimension1}
Let $M$ be a smooth connected manifold of dimension $m$ equipped
with a smooth distribution $\Delta$ of dimension $m-1$. If $M$
does not admit horizontal immersed submanifolds of dimension
$m-1$, then given two points $p,q$ in $M$, there exists at least
one piecewise smooth horizontal (with respect to $\Delta$) curve
connecting them (in fact $\Delta$ is equiregular).
\end{theorem}
\proof Let $\mathcal{F}$ be the set of all smooth vector fields
tangent to $\Delta$, that is,
$$
\mathcal{F}=\{f\in{\Gamma(TM)\mid f(p)\in{\Delta_{p}}}\textrm{ for
any } p\in{M}\}.
$$

For any $p\in{M}$ we claim that the orbit $\mathcal{O}_{p}$ of $p$
of the family $\mathcal{F}$ is of full dimension and thus
$\mathcal{O}_{p}$ is an open set of $M$.

In fact, if not, then $\textrm{dim}(\mathcal{O}_{p})\leq m-1$
where $\textrm{dim}(\mathcal{O}_{p})$ denotes the dimension of
$\mathcal{O}_{p}$. On the other hand by Corollary \ref{contained},
we have $\Delta_{q}\subset T_{q}\mathcal{O}_{p}$ for any
$q\in{\mathcal{O}_{p}}$ and so the dimension of $\mathcal{O}_{p}$
is not less than $\textrm{dim}(\Delta)=m-1$. Thus,
$\textrm{dim}(\mathcal{O}_{p})= m-1$ and
$T_{q}\mathcal{O}_{p}=\Delta_{q}$ for any $q\in{\mathcal{O}_{p}}$.
So $\mathcal{O}_{p}$ is a horizontal immersed submanifold of
dimension $m-1$. This contradicts with the assumption.

Since $M$ is connected and $M$ is the union of all orbits, we have
$M=\mathcal{O}_{p}$ for any $p\in{M}$.\endproof
\begin{lemma}\label{contactisometric}
Let $M$ be a smooth orientable contact manifold of dimension
$2n+1$ with a contact form $\eta$. Then $M$ does not admit
horizontal (with respect to $\Delta=\ker(\eta)$) immersion
submanifolds of dimension more than $n$.
\end{lemma}
\proof This is a well known fact. For the readers' convenience and
completeness we give a proof.

Let $i:N\rightarrow M$ be a horizontal immersed submanifold of $M$
and let $p\in{N}$. By definition,
$\dim(i_{*}(T_{p}N))=\dim(T_{p}N)$ and for any
$v_{1},v_{2}\in{i_{*}(T_{p}N)}$,
\begin{equation}\label{symplecticorthonogal}
\omega_{i(p)}(v_{1},v_{2})=0
\end{equation}
 where $\omega=d\eta$. Since
$\omega_{i(p)}$ is a symplectic form on the horizontal space
$\Delta_{p}$, from the non-degeneracy of $\omega_{i(p)}$ we have
that
$$
\dim((i_{*}(T_{p}N))^{\bot})+\dim(i_{*}(T_{p}N))=2n
$$
where $(i_{*}(T_{p}N))^{\bot}$ is  the symplectic orthogonal
subspace of $i_{*}(T_{p}N)$. \eqref{symplecticorthonogal} is
equivalent to $i_{*}(T_{p}N)\subset(i_{*}(T_{p}N))^{\bot}$. Thus
$$
2\dim{i_{*}(T_{p}N)}\leq
\dim((i_{*}(T_{p}N))^{\bot})+\dim(i_{*}(T_{p}N))=2n.
$$
\endproof
Now we are in the position to prove Theorem \ref{Main1}.
\proof[\textbf{Proof of Theorem \ref{Main1}}] Since $2n-1>n$, the
statement follows directly from Theorem \ref{codimension1} and Lemma
\ref{contactisometric} if $S$ is noncharacteristic. If $S$ possesses
characteristic points, let $\mathcal{F}^{\Delta}_{S}$ be the set of
all smooth vector fields tangent to $T^{\mathcal{H}}S$, that is,
$$
\mathcal{F}^{\Delta}_{S}:=\{f\in{\Gamma(TS)}\mid
f(p)\in{T_{p}^{\mathcal{H}}S}\textrm{ for any }p\in{S}\}.$$ Let
$\mathcal{O}_{p}$ be an orbit of $p\in{S}$ of the family
$\mathcal{F}^{\Delta}_{S}$. By Theorem \ref{orbit} $\mathcal{O}_{p}$
is an immersed submanifold. From the assumption that $\Sigma_{S}$
does not contain immersed submanifolds we conclude that there exists
at least one point $q\in{\mathcal{O}_{p}}$ which is not a
characteristic point. By Corollary \ref{contained} we have that
$\mathcal{F}^{\Delta}_{S}(q)\subset T_{q}(\mathcal{O}_{p})$ and the
dimension of $\mathcal{O}_{p}$ is not less than
$2n-1=\dim(\mathcal{F}^{\Delta}_{S}(q))$. Thus if
$\dim(\mathcal{O}_{p})=2n-1$, then $\mathcal{O}_{p}$ is a horizontal
immersion submanifold of $S$ (and $M$). This contradicts with Lemma
\ref{contactisometric} since $n>1$. So $\mathcal{O}_{p}$ is of full
dimension and it is an open set of $S$. The assertion follows from
the connectedness of $S$.
\endproof
\begin{example}[\textbf{the gauge ball in $\mathbb{H}^{n}$ with
$n>1$}] Let $\mathbb{S}^{n}=\{p\in{\mathbb{R}^{2n+1}}\mid
\|p\|=1\}$ where
$\|p\|:=\left((|x|^{2}+|y|^{2})^{2}+|t|^{2}\right)^{\frac{1}{4}}$
is the gauge norm in $\mathbb{H}^{n}$. $\mathbb{S}^{n}$ is called
the gauge ball centered in the origin. It is trivial to check that
the metric induced by the gauge norm is left-invariant and
1-homogeneous with respect to natural dilations in
$\mathbb{H}^{n}$. By direct computation the characteristic set
$\Sigma_{\mathbb{S}^{n}}$ of $\mathbb{S}^{n}$ consists of only two
points:
$$
\Sigma_{\mathbb{S}^{n}}=\{(0,0,1),(0,0,-1)\}.
$$
By Theorem \ref{Main1}, we see that the induced
Carnot-Carath\'{e}odory metric $\mathbb{S}^{n}$ can be defined.
\end{example}

\begin{example}[\textbf{hyperplanes in $\mathbb{H}^{n}$ with
$n>1$}] Let $n>1$. The vertical hyperplane
$L_{i}=\{(x_{1},\cdots,x_{i-1},0,x_{i+1},\cdots,x_{n},y,t)\in{R^{2n+1}}\}$
of $\mathbb{H}^{n}$ is a Lie subgroup of $\mathbb{H}^{n}$ with the
induced group law (that is the restriction to $L_{i}$ of the group
law of $\mathbb{H}^{n}$) and Lie algebra
$\mathcal{L}_{i}:=\overline{V}_{1}\bigoplus \overline{V}_{2}$
where
$\overline{V}_{1}=\textrm{span}\{X_{1},\cdots,X_{i-1},X_{i+1},\cdots,X_{n},
X_{n+1},\cdots,X_{n+i-1},\overline{X}_{n+i},X_{n+i+1},\cdots,X_{2n}\}$
and $\overline{V}_{2}=\textrm{span}\{T\}$ where
$\overline{X}_{n+i}:=\frac{\partial}{\partial y_{i}}$ and $X_{i},
T$ as in Example \ref{Heisenberg}). $L_{i}$ is noncharacteristic
and is a Carnot group. It is easy to prove that the
Carnot-Carath\'{e}odory metric $d_{\textrm{in}}$ induced by
$(\overline{V}_{1},g^{i}_{c})$ where $g^{i}_{c}$ is the
restriction of $g_{c}$ to $\overline{V}_{1}$ and the restricted
Carnot-Carath\'{e}odory metric $d_{\textrm{re}}$ satisfies that
$$
d_{\textrm{re}}\leq d_{\textrm{in}}\leq C d_{\textrm{re}}
$$
where $C$ is an absolute constant, see \cite{Gromov}, and there
are points $p,q\in{L_{i}}$ such that $d_{\textrm{re}}(p,q)<
d_{\textrm{in}}(p,q)$.

\end{example}
The following example shows the condition of dimension more than
three is unavoidable.
\begin{example}In the simplest Heisenberg group $\mathbb{H}^{1}$, we
consider horizontal curves in the horizontal hyperplane
$L_{t}=\{(x,y,0)\in{\mathbb{R}^{3}}\}$. Note that the point
$(0,0,0)$ is the unique characteristic point in $L_{t}$. Let
$\gamma(t)=(x(s),y(s),0)$ be a horizontal curve in $L_{t}$. Then
from \eqref{horizontal} we have
$$
x(s)\dot{y}(s)=y(s)\dot{x}(s)
$$
and hence
$$
y(s)=Cx(s)
$$
for some positive constant $C$. Thus there are no horizontal
curves in any hypersurface $S\subset L_{1}$ which does not contain
the point $(0,0,0)$. The same argument shows that there are no
horizontal curves in $L_{x}=\{(0,y,t)\in{\mathbb{R}^{3}}\}$ and in
$L_{y}=\{x,0,t)\in{\mathbb{R}^{3}}\}$.
\end{example}
\begin{remark}As in Euclidean geometry, $\mathbb{R}^{3}$ can serve as a model to the study of higher
dimension, in the study of sub-Riemannian geometry
$\mathbb{H}^{1}$ can also be seen as a model. But in the
developing geometric measure  theory in the setting of
sub-Riemannian geometry, in particular for the notion of
rectifiability ( and possibly for co-area formulae) in Carnot
groups, $\mathbb{H}^{1}$ may be an exception. We recall that in
\cite{FSSC2}-\cite{FSSC5} B. Franchi, R. Serapioni and F. Serra
Cassano have proposed a notion of rectifiability by introducing
the notion of intrinsic regular hypersurfaces, and another notion
which is a counterpart of Federer's definition of rectifiability
where the ``model spaces'' are replaced by Carnot groups is
announced by S. Pauls in \cite{Pa1}. It seems that the Pauls'
notion is meaningless for $\mathbb{H}^{1}$. One reason is that the
codimension one Lie subgroups $L_{x},L_{y}$ of $\mathbb{H}^{1}$
have no stratified structure. Another reason is the horizontal
non-connectivity of hypersurfaces in $\mathbb{H}^{1}$ as shown in
the last example. A very intriguing question arises: whether are
the two notions of rectifiability for $\mathbb{H}^{n}(n>1)$
equivalent in any reasonable sense?
\end{remark}
\section{Horizontal connection,
the horizontal mean curvature and the horizontal divergence
theorem}\label{sec:proof2}

\begin{definition}[\textbf{horizontal connection}]
Let $g$ be any orthogonal extension of $g_{c}$ and let $\nabla$ be
the Levi-Civita connection with respect to $g$. We define the
horizontal connection $\widetilde{D}$ on $\Delta$ as
$$
\begin{aligned}
&\widetilde{D}:\Gamma(\Delta)\times\Gamma(\Delta)\rightarrow \Gamma(\Delta)\\
&\widetilde{D}_{X}Y=\sum_{i=1}^{k}g(\nabla_{X}Y,X_{i})X_{i}\textrm{
for any } X,Y\in{\Gamma(\Delta)}
\end{aligned}
$$
where $\{X_{1},\cdots,X_{k}\}$ is an orthonormal basis of
$\Delta$.
\end{definition}
\begin{remark}\label{independence}The definition of $\widetilde{D}$ is independent of the choice
of orthonormal basis of $\Delta$. In fact, let
$\overline{X}_{i}=\sum_{j=1}^{k}a_{ij}X_{j},i=1,\cdots,k$ be
another orthonormal basis. Then $(a_{ij})$ is an orthonormal
matrix (everywhere) and hence
$$
\begin{aligned}
\sum_{i=1}^{k}g(\nabla_{X}Y,\overline{X}_{i})\overline{X}_{i}
&=\sum_{j=1}^{k}\sum_{l=1}^{k}\sum_{i=1}^{k}a_{ij}a_{il}g(\nabla_{X}Y,X_{j})X_{l}
=\sum_{j=1}^{k}\sum_{l=1}^{k}\delta_{jl}g(\nabla_{X}Y,X_{j})X_{l}\\
&=\sum_{j=1}^{k}g(\nabla_{X}Y,X_{j})X_{j}
\end{aligned}
$$
Thus $\widetilde{D}$ is well defined.
\end{remark}

We call $\widetilde{D}$ a ``connection'' because of the following
fact.
\begin{lemma}\label{fundamentalpro}
$\widetilde{D}$  satisfies the following properties
\begin{enumerate}
    \item\label{1} $\widetilde{D}_{X}Y$ is $R-$linear in both arguments,
    \item\label{2} $\widetilde{D}_{X}Y$ is $C^{\infty}(M)-$linear in the argument of
    $X$,
    \item\label{4} the Leibniz rule holds:
    $$
    \widetilde{D}_{X}(fY)=(Xf)Y+f\widetilde{D}_{X}Y\textrm{ for any }
    f\in{C^{\infty}(M)}, X,Y\in{\Gamma(\Delta)}
    $$
    \item\label{4''} $\widetilde{D}$ is compatible with respect to $g_{c}$ , that
    is,
   \begin{equation}\label{compability}
    Xg_{c}(Y,Z)=g_{c}(\widetilde{D}_{X}Y,Z)+g_{c}(Y,\widetilde{D}_{X}Z)\textrm{ for any
    }X,Y,Z\in{\Gamma(\Delta)},
   \end{equation}
   \item\label{5} the following ``symmetry'' property holds:
     \begin{equation}\label{symmetry}
 \widetilde{D}_{X}Y-\widetilde{D}_{Y}X=[X,Y]^{\mathcal{H}}\textrm{ for any }X,Y\in{\Gamma(\Delta)}
\end{equation}
   where $[X,Y]^{\mathcal{H}}$ is, by definition, the projection
   of $[X,Y]$ to $\Delta$, that is,
   $$
   [X,Y]^{\mathcal{H}}=\sum_{i=1}^{k}g([X,Y],X_{i})X_{i}.
   $$
\end{enumerate}
\end{lemma}
\proof The proof is trivial. All follow directly from the
definition of $\widetilde{D}$, the compatibility and the symmetry
of the Levi-Civita connection $\nabla$ together with the fact that
$g$ is an orthogonal extension of $g_{c}$.
\endproof
\begin{remark}\label{projection}For any vector field $X$, the horizontal part $X^{\mathcal{H}}$
of $X$: $X^{\mathcal{H}}=\sum_{i=1}^{k}g(X,X_{i})X_{i}$, is
independent of  any orthogonal extension $g$ of $g_{c}$. In fact,
since $g$ is an orthogonal extension of $g_{c}$, the projection,
with respect to the decomposition of direct summation, of a vector
field to the horizontal bundle is the same as the projection of
this vector field to the horizontal bundle, with respect to the
orthogonal decomposition.

From \eqref{2} and \eqref{4} of Lemma \ref{fundamentalpro} it is
straight to verify that $\widetilde{D}_{X}Y(p)$ depends only on
$X(p)$ and the evaluations of $Y$ in a neighborhood of $p$.

Any operator from
$\Gamma(\Delta)\bigotimes\Gamma(\Delta)\rightarrow \Gamma(\Delta)$
satisfying \eqref{1}, \eqref{2}, \eqref{4} and \eqref{4''} of
Lemma \ref{fundamentalpro} is called a \textit{nonholonomic
connection on $\Delta$}.
\end{remark}
\begin{proposition}\label{independet}
The operator $\widetilde{D}$ is independent of the choice of
orthogonal extensions of $g_{c}$.
\end{proposition}
\proof From Lemma \ref{fundamentalpro} we know that
$\widetilde{D}_{U}V=\sum_{j=1}^{k}U(V^{j})X_{j}+\sum_{i=1}^{k}\sum_{j=1}^{k}U^{i}V^{j}\widetilde{D}_{X_{i}}X_{j}$
for any $U=\sum_{i=1}^{k}U^{i}X_{i},\,
V=\sum_{j=1}^{k}V^{j}X_{j}\in{\Gamma(\Delta)}$. Thus
$\widetilde{D}$ is determined by the connection coefficients
$\Gamma_{ij}^{l}=g_{c}(\widetilde{D}_{X_{i}}X_{j},X_{l})=g(\nabla_{X_{i}}X_{j},X_{l})$.
By the Cozhul's formula, for any $i,j,l=1,\cdots,k$ we have
$$
\begin{aligned}
g(\nabla_{X_{i}}X_{j},X_{l})=&\frac{1}{2}\{X_{i}g(X_{j},X_{l})+X_{j}g(X_{l},X_{i})-X_{l}g(X_{i},X_{j})\\
&\quad
-g(X_{i},[X_{l},X_{j}])-g(X_{j},[X_{i},X_{l}])+g(X_{l},[X_{j},X_{i}])\}\\
=&-\frac{1}{2}\{g(X_{i},[X_{l},X_{j}])+g(X_{j},[X_{i},X_{l}])-g(X_{l},[X_{j},X_{i}])\}\\
=&-\frac{1}{2}\{g_{c}(X_{i},[X_{l},X_{j}]^{\mathcal{H}})+
g_{c}(X_{j},[X_{i},X_{l}]^{\mathcal{H}})-g_{c}(X_{l},[X_{j},X_{i}]^{\mathcal{H}})\},
\end{aligned}
$$
since $\{X_{1},\cdots,X_{k}\}$ is an orthonormal basis of $\Delta$
with respect to  $g_{c}$. The assertion follows from Remark
\ref{projection}. \endproof \vskip10pt

Now we turn to the proof of Theorem \ref{Main2}. The definition of
$D$ in Theorem \ref{Main2} is inspired by the proof of Proposition
\ref{independet} and the fundamental theorem in Riemannian
geometry.

\proof[\textbf{Proof of Theorem \ref{Main2}}]In this proof, to
simplify the notations we use the Einstein  summation convention:
if in any term the same index name appears twice, as both an upper
and a lower index, that term is assumed to be summed over all
possible values of that index (from 1 to $k$). We first prove that
$D$ is independent of the choice of orthonormal basis of $\Delta$.
Let $\{\widetilde{X}_{1},\cdots,\widetilde{X}_{k}\}$ be another
orthonormal basis of $\Delta$, where
$\widetilde{X}_{i}=a_{i}^{j}X_{j}$, $i=1,\cdots,k$ and
$A=(a_{i}^{j})$ is an orthogonal matrix (everywhere). For any
$U=\widetilde{U}^{j}\widetilde{X}_{j}=U^{j}X_{j},$
$V=\widetilde{V}^{i}\widetilde{X}_{i}=V^{i}X_{i}\in{\Gamma(\Delta)}$
where $\widetilde{U}^{j}=b_{r}^{j}U^{r}$,
$\widetilde{V}^{i}=b_{s}^{i}V^{s}$ and $B=(b_{i}^{j})$ is the
inverse matrix of $A$ with $b_{i}^{j}=a_{j}^{i}$  for
$i,j=1,\cdots,k$ since $A$ is orthogonal, we try to compute
$$
U(\widetilde{V}^{i})\widetilde{X}_{i}+
\widetilde{U}^{j}\widetilde{V}^{i}
\widetilde{\Gamma}_{ij}^{l}\widetilde{X}_{l}=\textrm{I}+\textrm{II}
$$
where $
\widetilde{\Gamma}_{ij}^{l}=-\frac{1}{2}\{g_{c}(\widetilde{X}_{i},[\widetilde{X}_{j},\widetilde{X}_{l}]^{\mathcal{H}})+
g_{c}(\widetilde{X}_{j},[\widetilde{X}_{l},\widetilde{X}_{i}]^{\mathcal{H}})
-g_{c}(\widetilde{X}_{l},[\widetilde{X}_{i},\widetilde{X}_{j}]^{\mathcal{H}})\}.
$ Since
$$
[\widetilde{X}_{j},\widetilde{X}_{l}]=a_{j}^{r}a_{l}^{h}
[X_{r},X_{h}]+a_{j}^{r}X_{r}(a_{l}^{h})X_{h}-a_{l}^{h}X_{h}(a_{j}^{r})X_{r},
$$
we have
$$
[\widetilde{X}_{j},\widetilde{X}_{l}]^{\mathcal{H}}=a_{j}^{r}a_{l}^{h}
[X_{r},X_{h}]^{\mathcal{H}}+a_{j}^{r}X_{r}(a_{l}^{h})X_{h}-a_{l}^{h}X_{h}(a_{j}^{r})X_{r}
$$
and hence
$$
g_{c}(\widetilde{X}_{i},[\widetilde{X}_{j},\widetilde{X}_{l}]^{\mathcal{H}})
=a_{i}^{s}a_{j}^{r}a_{l}^{h}g_{c}(X_{s},[X_{r},X_{h}]^{\mathcal{H}})+
\sum_{s=1}^{k}a_{i}^{s}a_{j}^{r}X_{r}(a_{l}^{s})-\sum_{s=1}^{k}a_{i}^{s}a_{l}^{h}X_{h}(a_{j}^{s}),
$$
using the same arguments to other terms, and from
$b_{i}^{j}a_{l}^{i}=\delta_{l}^{j}$,
$a_{l}^{i}a_{l}^{j}=\delta_{i}^{j}$ for $i,j=1,\cdots,k$, we
deduce
$$
\begin{aligned}
\textrm{II}=&-\frac{1}{2}\sum_{h=1}^{k}V^{s}U^{r}\{g_{c}(X_{s},[X_{r},X_{h}]^{\mathcal{H}})
+g_{c}(X_{r},[X_{h},X_{s}]^{\mathcal{H}})-g_{c}(X_{h},[X_{s},X_{r}]^{\mathcal{H}})\}X_{h}\\
&-\frac{1}{2}\{\sum_{s=1}^{k}\sum_{l=1}^{k}U^{r}V^{s}a_{l}^{h}X_{r}(a_{l}^{s})X_{h}-\sum_{s=1}^{k}
\sum_{h=1}^{k}U^{r}V^{s}b_{r}^{j}X_{h}(a_{j}^{s})X_{h}-
U^{r}V^{s}b_{r}^{j}X_{s}(a_{j}^{h})X_{h}
\\&\qquad\,+U^{r}V^{s}b_{s}^{i}X_{r}(a_{i}^{h})X_{h}
-\sum_{r=1}^{k}\sum_{h=1}^{k}U^{r}V^{s}b_{s}^{i}X_{h}(a_{i}^{r})X_{h}+
\sum_{r=1}^{k}\sum_{l=1}^{k}U^{r}V^{s}a_{l}^{h}X_{s}(a_{l}^{r})X_{h}\}.
\end{aligned}
$$
Renaming and rearranging some indices, then cancelling or adding
some terms (using $X_{l}(a_{r}^{j}b_{j}^{h})=0$ and
$a_{l}^{r}=b_{r}^{l}$ for any $l,r,h=1,\cdots,k$), we infer that
$$
\textrm{II}=U^{r}V^{s}\Gamma_{sr}^{h}X_{h}-U^{r}V^{j}a_{i}^{l}X_{r}(b_{j}^{i})X_{l}.
$$
Since
$$
\begin{aligned}
\textrm{I}&=U(\widetilde{V}^{j})\widetilde{X}_{j}=U^{r}X_{r}(b_{i}^{j}V^{i})a_{j}^{l}X_{l}\\
&=U^{r}X_{r}(V^{i})b_{i}^{j}a_{j}^{l}X_{l}+U^{r}V^{i}X_{r}(b_{i}^{j})X_{l}\\
&=U(V^{j})X_{j}+U^{r}V^{i}X_{r}(b_{i}^{j})a_{j}^{l}X_{l},
\end{aligned}
$$
we obtain
$$
U(\widetilde{V}^{i})\widetilde{X}_{i}+
\widetilde{U}^{j}\widetilde{V}^{i}
\widetilde{\Gamma}_{ij}^{l}\widetilde{X}_{l}=U(V^{j})X_{j}+U^{r}V^{s}\Gamma_{sr}^{h}X_{h}.
$$
Thus $D$ is well defined.

From the definition of $D$ we have for any $i,j,l=1,\cdots,k$,
\begin{align}
&\Gamma_{ij}^{l}=g_{c}(D_{X_{j}}X_{i},X_{l}),\label{equal}\\
&\Gamma_{ji}^{l}+\Gamma_{li}^{j}=0,\label{zero}\\
&\Gamma_{ij}^{l}-\Gamma_{ji}^{l}=g_{c}([X_{j},X_{i}]^{\mathcal{H}},X_{l}).\label{antisymmetric}
\end{align}

 Properties
\eqref{1'}-\eqref{3'} are obvious. To prove \eqref{4'}, let
$U=U^{j}X_{j},V=V^{i}X_{i},W=W^{l}X_{l}$. Then from \eqref{2'} and
\eqref{3'}  we have
$$
\begin{aligned}
g_{c}(D_{V}U,W)+g_{c}(U,D_{V}W)&=U^{j}W^{l}g_{c}(D_{X_{j}}(V^{i}X_{i}),X_{l})
+U^{j}V^{i}g_{c}(X_{j},D_{X_{i}}(W^{l}X_{l}))\\
&=U^{j}V^{i}W^{l}(g_{c}(D_{X_{j}}X_{i},X_{l})+g_{c}(X_{j},D_{X_{i}}X_{l}))+\sum_{i=1}^{k}W^{i}V(U^{i})\\
&\quad+\sum_{j=1}^{k}U^{j}V(W^{j})\\
&=U^{j}V^{i}W^{l}(g_{c}(D_{X_{j}}X_{i},X_{l})+g_{c}(X_{j},D_{X_{i}}X_{l}))+Vg_{c}(W,U).
\end{aligned}
$$
Since by \eqref{equal} and \eqref{zero}
$$
\begin{aligned}
g_{c}(D_{X_{j}}X_{i},X_{l})+g_{c}(X_{j},D_{X_{i}}X_{l})&=\Gamma_{ji}^{l}+\Gamma_{lj}^{i}\\
&=0
\end{aligned}
$$
for any $i,j,l=1,\cdots,k$, we conclude that
$$
g_{c}(D_{V}U,W)+g_{c}(U,D_{V}W)=Vg_{c}(W,U).
$$

To see \eqref{symmetry2}, for any $l=1,\cdots,k$, $U=U^{j}X_{j},
V=V^{i}X_{i}\in{\Gamma(\Delta)}$, by \eqref{antisymmetric} we have
$$
\begin{aligned}
g_{c}(D_{U}V-D_{V}U,X_{l})&=U(V^{l})-V(U^{l})+U^{j}V^{i}(\Gamma_{ij}^{l}-\Gamma_{ji}^{l})\\
&=U(V^{l})-V(U^{l})+U^{j}V^{i}g_{c}([X_{j},X_{i}]^{\mathcal{H}},X_{l})\\
&=g_{c}([U,V]^{\mathcal{H}},X_{l}).
\end{aligned}
$$

The uniqueness follows from an argument similar to the Riemannian
case: first use \eqref{2'} and \eqref{3'} to obtain
\eqref{definition} where
$\Gamma_{ij}^{l}=g_{c}(D_{X_{j}}X_{i},X_{l})$, then using
\eqref{4'} and \eqref{symmetry2} to get a formula of Cozhul type,
$g_{c}(D_{X_{j}}X_{i},X_{l})=-\frac{1}{2}\{g_{c}(X_{i},[X_{l},X_{j}]^{\mathcal{H}})+
g_{c}(X_{j},[X_{i},X_{l}]^{\mathcal{H}})-g_{c}(X_{l},[X_{j},X_{i}]^{\mathcal{H}})\}.
$
\endproof
\begin{remark} Thus $D$ is an intrinsic notion of the
sub-Riemannian structure $(\Delta,g_{c})$. Here we must point out
that the connection $\overline{D}$ obtained by projecting the
Levi-Civita connection onto the horizontal bundle is first
introduced by Schouten in \cite{Sch}\footnote{There the Riemannian
metric $g$ with respect to which the Levi-Civita connection was
computed is not necessary to be an orthogonal extension of
$g_{c}$. Thus in this case $\overline{D}$ is not intrinsic: it
depends on $g$.}, further developed by Vagner in
\cite{Vagner1}-\cite{Vagner4} (see \cite{Gor} for a survey), see
also \cite{VF} and \cite{VG} where $\overline{D}$ is called
``truncated connection''. But it seems that they did not realize
the existence of $D$ in Theorem \ref{Main2} and its intrinsic
nature (determination by $(\Delta,g_{c})$) of $D$, see p.202 in
\cite{VG}.

In Russian mathematicians' papers (also E. Cartan's, see
\cite{Cartan}, \cite{BKMM} and \cite{KRP}) nonholonomic connection
was mainly used to study ``geodesics'' (but not sub-Riemannian
geodesics i.e. nonholonomic geodesics) and curvature of
distributions in the setting of nonholonomic dynamic systems. Let
$\gamma(s)$ be a smooth horizontal curve in $M$. We call $\gamma$
a ``geodesic'' if $D_{\dot{\gamma}}\dot{\gamma}=0$. The mechanical
significance of such ``geodesics'' is that they characterize the
trajectories of motion of a mechanical system with quadratic
Lagrangian and linear constraints (say $\Delta$), see \cite{VF}
and \cite{VG} for details.
\end{remark}
\begin{example}[\textbf{the horizontal connection in Carnot
group}]\label{connectioncarnot} In a Carnot group $G$, since its
Lie algebra $\mathcal{G}$ is graded, if we choose a system of left
invariant vector fields $\{X_{1},\cdots,X_{k}\}$ as an orthogonal
basis of the horizontal bundle $\Delta=V_{1}$, then the connection
coefficients $\Gamma_{ij}^{l}$ of the horizontal connection $D$
are vanishing everywhere for any $i,j,l=1,\cdots,k$. For this
basis, $D$ has the simple form
$$
D_{U}V=\sum_{i=1}^{k}U(V^{i})X_{i}\textrm{ for any }U,
V=\sum_{i=1}^{k}V^{i}X_{i}\in{\Gamma(\Delta)},
$$
see \cite{TY} for details.
\end{example}

For completeness we introduce the notion of horizontal divergence.
\begin{definition}[\textbf{horizontal divergence}]Let
$X\in{\Gamma(\Delta)}$. The horizontal divergence
$\textrm{div}_{\mathcal{H}}X$ of $X$ is defined by
$$
\textrm{div}_{\mathcal{H}}X=\sum_{i=1}^{k}g_{c}(D_{X_{i}}X,X_{i})
$$
Note that as shown in Remark \ref{independence},
$\textrm{div}_{\mathcal{H}}X$ is independent of the choice of
orthonormal basis of $\Delta$ because of \eqref{2} in Lemma
\ref{fundamentalpro}.
\end{definition}
The following proposition, which in the case of Carnot groups are
well known, follows immediately from the definitions.
\begin{proposition}\label{divergence}
$\textrm{div}_{\mathcal{H}}X=\textrm{div}X$ for any
$X\in{\Gamma(\Delta)}$ where $\textrm{div}X$ is the usual
divergence of $X$ computed with respect to any orthogonal
extension $g$ of $g_{c}$.
\end{proposition}
\vskip10pt

 It is interesting that $D$ can be used to define the
horizontal mean curvature of (noncharacteristic) hypersurfaces in
$(M,\Delta,g_{c})$. In the rest of this paper, we discuss how this
can be done. Roughly speaking, in sub-Riemannian (or nonholonomic)
geometry horizontal connection, horizontal tangent connection,
horizontal normal, horizontal second fundamental form and
horizontal mean curvature are counterparts of Levi-Civita
connection, tangent connection, Riemannian normal, second
fundamental form and mean curvature in Riemannian geometry,
respectively. \vskip10pt

We assume, if without further notice, $S$ is a smooth
\textbf{noncharacteristic} hypersurface in a sub-Riemannian
manifold $(M,\Delta,g_{c})$. Then $T^{\mathcal{H}}S$ is a
subbundle of $TS$ of dimension $k-1$. From the definition we have
\begin{lemma}\label{selfcontained}
If $X,Y\in{\Gamma(T^{\mathcal{H}}S)}$, then
$[X,Y]^{\mathcal{H}}\in{\Gamma(T^{\mathcal{H}}S)}$.
\end{lemma}

It is clear that any vector $v$ in $T_{p}^{\mathcal{H}}S
(p\in{S})$ can be extended to a vector field in $T^{\mathcal{H}}S$
by first extending $v$ to a vector field $V$ in $TS$ then
projecting $V$ to $T^{\mathcal{H}}S$, and any  vector field $V$ in
$T^{\mathcal{H}}S$ can smoothly extended to a horizontal vector
field in $\Delta$ by first extending $V$ to a vector field
$\overline{V}$ in $TM$ then projecting $\overline{V}$ to $\Delta$.
Sometimes we will denote by the same symbol both the extended
vector field and the original vector field, in particular when we
have confirmed that the objects under consideration are
independent of such extensions.

If $X,Y$ are vector fields in $T^{\mathcal{H}}S$, we can extend
them to horizontal vector fields $\overline{X},\overline{Y}$ in
$M$, apply the ambient derivative operator $D$, and then decompose
at points of $S$ to get
\begin{equation}\label{eq:decomposition}
D_{\overline{X}}\overline{Y}(x)=(D_{\overline{X}}\overline{Y})^{\top}(x)+(D_{\overline{X}}\overline{Y})^{\bot}(x)
\quad x\in{S}
\end{equation}
where $(D_{\overline{X}}\overline{Y})^{\top}(x)$,
$(D_{\overline{X}}\overline{Y})^{\bot}(x)$ are the projections of
$D_{\overline{X}}\overline{Y}(x)$ onto $T^{\mathcal{H}}_{x}S$ and
the direction of $\mathcal{V}(x)$ respectively, where
$\mathcal{V}(x)$ is the unit horizonal normal, see Remark
\ref{normalintrinsic}.

\begin{definition}[\textbf{the horizontal second fundamental form}]
Let $X,Y$ be vector fields in $T^{\mathcal{H}}S$. We define
$$
\textrm{II}(X,Y):=(D_{\overline{X}}\overline{Y})^{\bot}
$$
where $\overline{X},\overline{Y}$ are the arbitrarily extended
horizontal vector fields of $X,Y$ respectively, as \textit{the
horizontal second fundamental form} of $S$.
\end{definition}

\begin{theorem}\label{symmetric}
The horizontal second fundamental form $\textrm{II}(X,Y)$ is
    \begin{enumerate}
    \item independent of the extension of $X$ and $Y$;
    \item bilinear over $C^{\infty}(S)$; and
    \item\label{3} symmetric in $X$ and $Y$.
    \end{enumerate}
\end{theorem}
The proof of Theorem \ref{symmetric} is very similar to the
Riemannian case, see \cite{TY}. Note that \eqref{3} follows from
Lemma \ref{selfcontained} and the symmetry property
\eqref{symmetry} of $D$, as in the Riemannian geometry the
symmetry of the second fundamental follows directly from the
symmetry (torsion free) of the Levi-Civita connection. Theorem
\ref{symmetric} in particular implies that $\textrm{II}(X,Y)(p)$
depends only on $X(p)$ and $Y(p)$.
\begin{definition}[\textbf{the horizontal tangent connection}] We
define the horizontal tangent connection
$$
D^{\top}:\Gamma(T^{\mathcal{H}}S)\bigotimes\Gamma(T^{\mathcal{H}}S)\rightarrow
\Gamma(T^{\mathcal{H}}S)
$$
by
$$
D^{\top}_{X}Y=(D_{\overline{X}}\overline{Y})^{\top},
$$
where $\overline{X},\overline{Y}$ are the arbitrarily extended
horizontal vector fields of $X,Y$ respectively.
\end{definition}

\begin{theorem}\label{tangentconnection}
The horizontal tangent connection $D^{\top}$ is well defined, that
is, $D^{\top}_{X}Y$ is independent of the extension of $X$ and
$Y$. Moreover $D^{\top}$ satisfies \eqref{1}-\eqref{5} of Lemma
\ref{fundamentalpro} where $\Delta$ is replaced by
$T^{\mathcal{H}}S$ and $g_{c}$ is replaced by
$g_{c}^{\mathcal{H}}$ which is the restriction of $g_{c}$ to
$T^{\mathcal{H}}S$.
\end{theorem}
The proof of Theorem \ref{tangentconnection} is direct and
trivial. From the preceding discussions, we know that in $S$ there
exists an intrinsic horizontal connection determined by the
sub-Riemannian structure $(T^{\mathcal{H}}S,g_{c}^{\mathcal{H}})$,
since
$$
TS=T^{\mathcal{H}}S\stackrel{\overline{g}}\bigoplus
(TS/T^{\mathcal{H}}S)
$$
where $\overline{g}$ is the restriction of $g$ to $TS$. From the
definitions it is easy to see that $\overline{g}$ is the orthogonal
extension of $g^{\mathcal{H}}_{c}$, even if $T^{\mathcal{H}}S$ may
not satisfy the Chow condition. Thus the ``external'' connection
$D^{\top}$ is equal to the intrinsic one by Theorem
\ref{tangentconnection} and Theorem \ref{Main2}.
\begin{definition}[\textbf{horizontal tangent divergence}] Let
$X\in{\Gamma(T^{\mathcal{H}}S)}$. The horizontal tangent
divergence $div_{\mathcal{H}}^{\top}X$ of $X$ is defined by
$$
div_{\mathcal{H}}^{\top}X=\sum_{i=1}^{k-1}g_{c}^{\mathcal{H}}(D^{\top}_{\tau_{i}}X,\tau_{i})
$$
where $\{\tau_{1},\cdots,\tau_{k-1}\}$ is a orthonormal basis of
$T^{\mathcal{H}}S$ with respect to $g_{c}^{\mathcal{H}}$. For a
horizontal vector field $Y$ defined in $S$, that is,
$Y\in{\Gamma(\Delta|_{S})}$ where $\Delta|_{S}$ is the restriction
of $\Delta$ to $S$, we also can define its horizontal tangent
divergence
$$
div_{\mathcal{H}}^{\top}Y=\sum_{i=1}^{k-1}g_{c}(D_{\tau_{i}}X,\tau_{i}).
$$
\end{definition}
Similar to Proposition \ref{divergence}, we have
\begin{proposition}Let $Y\in{\Gamma(T^{\mathcal{H}}S)}$. Then
$$
div_{\mathcal{H}}^{\top}Y=\textrm{div}_{S}Y
$$
where $\textrm{div}_{S}Y$ is the Riemannian tangent divergence of
$Y$ on $S$ computed with respect to any orthogonal extension $g$
of $g_{c}$.
\end{proposition}
\begin{definition}
The \textit{scalar horizontal second fundamental form} $h$ is the
symmetric bilinear function on $T^{\mathcal{H}}S$ defined by
$$
h(X,Y)=g_{c}(\textrm{II}(X,Y),\mathcal{V}).
$$
That is  $\textrm{II}(X,Y)=h(X,Y)\mathcal{V}$. Thus $h$ uniquely
determines an endomorphism of $T^{\mathcal{H}}S$, say $A$, that
is,
$$
g_{c}(AX,Y)=h(X,Y)\quad\textrm{ for all
}X,Y\in{\Gamma(T^{\mathcal{H}}S)}.
$$
$A$ is self-adjoint and we call $A$ \textit{the horizontal shape
operator} of $S$.
\end{definition}\vskip10pt
For any $p\in{S}$, $A$ gives a symmetric linear map
$A_{p}:T^{\mathcal{H}}_{p}S\rightarrow T^{\mathcal{H}}_{p}S$. Then
by the symmetry of $A_{p}$, $A_{p}$ has $k-1$ real eigenvalues.
\begin{definition}[\textbf{the horizontal mean curvature}]\label{horizontal meancurvature}
The $k-1$ eigenvalues of $A_{p}$,
$\kappa_{1},\cdots,\kappa_{k-1}$, are called \textit{the
horizontal principal curvatures} at $p$ and the corresponding
eigenspaces are called \textit{horizontal principal directions}.
We define the \textit{horizontal mean curvature} $H_{X}(p)$ at $p$
the trace of $A_{p}$, that is, $
H_{X}(p)=\sum_{i=1}^{k-1}\kappa_{i}$ and call the product of
$\kappa_{1},\cdots,\kappa_{k-1}$ the \textit{horizontal Gaussian
curvature} at $p$.
\end{definition}
It has been proven in \cite{TY} that in the case of Carnot groups
our definition of the horizontal mean curvature coincides with
that in \cite{DGN1} and \cite{Pa}.
\begin{remark}Since $D$, $T^{\mathcal{H}}S$ and $\mathcal{V}$ are intrinsic, so
is the notion of horizontal mean curvature. As in \cite{DGN1} and
\cite{Pa}, the horizontal mean curvature can be defined only at
noncharacteristic points\footnote{In \cite{DGN1} (see also
\cite{GP}), if $p$ is a characteristic point, the horizontal mean
curvature at $p$ is defined by $H(p)=\lim_{q\rightarrow p}H(q)$ if
this limit exists. Since it is impossible to determine the size of
the set of characteristic points at which the limits exist, this
definition seems meaningless}. However by Proposition
\ref{smallness}, the horizontal mean curvature can be defined
almost everywhere.
\end{remark}
\begin{example}
As in Example \ref{connectioncarnot}, choose a system of
left-invariant vector fields $\{X_{1},\cdots,X_{k}\}$ as an
orthogonal basis of the horizontal bundle, then for any
noncharacteristic hypersurface $S\subset G$, the horizontal mean
curvature can be expressed as the form
$$
H(p)=\sum_{i=1}^{k}X_{i}(\mathcal{V}^{i})(p)\textrm{ for any
}p\in{S}
$$
where $\mathcal{V}=\sum_{i=1}^{k}\mathcal{V}^{i}X_{i}$ is the unit
horizontal normal of $S$, see \cite{TY}.
\end{example}


\begin{thebibliography}{ABC}
\bibitem{Be}A. Bella\"{i}che, Tangent spaces in sub-Riemannian
geometry, in \textit{Sub-Riemannian Geometry}, 1-78, Prog. in
Math., v.144, edited by A. Bella\"{i}che and J. J. Risler,
Birkh\"{a}user, 1996.
\bibitem{BKMM}Anthony M. Bloch, P. S. Krishnaprasad, Jerrold E. Marsden, Richard M.
Murray, Nonholonomic mechanical systems with symmetry, Arch.
Rational Mech. Anal. \textbf{136}(1996), 21-99.
\bibitem{Cartan} E. Cartan. Sur la repres¡äentation g¡äeom¡äetrique des syst`emes
mat¡äeriels non holonomes, Proc. Int. Congr. Math., vol.
\textbf{4}(1928), Bologna, 253-261.
\bibitem{Chow} W. L. Chow, \"{u}ber Systeme non linearen partiellen
Differentialgleichungen erster Ordung. Math. Ann.
\textbf{117}(1939), 98-105.
\bibitem{DGN}D. Danielli, N. Garofalo, D. M. Nhien, Minimal
surfaces, surfaces of constant mean curvature and isoperimetry in
Carnot groups, 2001, preprint.
\bibitem{DGN1} D. Danielli, N. Garofalo, D. M. Nhien, Non-doubling
Ahlfors measures, perimeter measures, and the characterization of
the trace spaces of Sobolev functions in Carnot-Carath\'{e}odory
spaces. Memo. AMS. to appear.
\bibitem{De1}M. Derridj, Un probl\'{e}me aux limites pour une
classe d'op\'{e}rateurs du second ordre hypoelliptiques, Ann.
Inst. Fourier, Grenoble, \textbf{21}(1971), 99-148.
\bibitem{De2}M. Derridj, Sur un th\'{e}or\`{e}me de traces, Ann.
Inst. Fourier, Grenoble, \textbf{22}(1972), 73-83.
\bibitem{FSSC1}B. Franchi, R. Serapioni, F. Serra Cassano,
 Meyers-Serrin type theorems and relaxation of variational integrals
  depending on vector fields. Houston Math. J. \textbf{22}(1996), 858-889.
\bibitem{FSSC2}B. Franchi, R. Serapioni, F. Serra Cassano,
Rectifiability and perimeter in the Heisenberg group, Math. Ann.,
\textbf{321}(2001), 479-531.
\bibitem{FSSC3}B. Franchi, R. Serapioni, F. Serra Cassano, Regular
hypersurfaces, intrinsic perimeter and implicit functions theorem
in Carnot groups, to appear on Comm. in Analysis and
Geometry(2001), available from http://cvgmt.sns.it.
\bibitem{FSSC4}B. Franchi, R. Serapioni, F. Serra Cassano, On the
structure of finite perimeter sets in step 2 Carnot groups, 2001,
preprint, available from http://cvgmt.sns.it.
\bibitem{FSSC5}B. Franchi, R. Serapioni, F. Serra Cassano,
Rectifiability and perimeter in step 2 groups, 2001, preprint,
available from http://cvgmt.sns.it.
\bibitem{GN}N. Garofalo, D. M. Nhieu, Isoperimetric and Sobolev
inequalities for Carnot-Carath\'{e}odory spaces and the existence
of minimimal surfaces, Comm. Pure Appl. Math. \textbf{49}(1996),
1081-1144.
\bibitem{GP}N. Garofalo, S. D. Pauls, The Bernstein problem in the
Heisenberg group, 2002, preprint, available from:
http://cvgmt.sns.it.
\bibitem{Gor} E. M. Gorbatenko, Differential geometry of
nonholonomic manifolds (after V.V.Vagner), Geom. Sbornik, 31-43,
Tomsk University, Tomsk, 1985.
\bibitem{Gromov}M. Gromov, Carnot-Carath\'{e}odory spaces seen from
within, in \textit{Sub-Riemannian Geometry}, 79-323, Prog. in
Math., v.144, edited by A. Bella\"{i}che and J. J. Risler,
Birkh\"{a}user, 1996.
\bibitem{KRP}J. Koiller, P. R. Rodrigues, P. Pitanga, Non-holonomic connections
following \'{E}lie Cartan, An. Acad. bras. Cienc. \textbf{7}(2001)
165-190.
\bibitem{Ko}U. Koschorke, \textit{Vector fields and other vector bundle
morphisms-a singularity approach}. Lecture notes in Mathematics
847, Springer, Berlin, 1981.
\bibitem{Magnani}V. Magnani, \textit{Elements of geometric measure theory
on sub-Riemannian groups}, Ph.D. thesis, Scuola Normale Superiore,
Pisa, 2002.
\bibitem{Magnani1}V. Magnani, A blow-up theorem for regular hypersurfaces on nilpotent groups, Manuscripta Math.,
 \textbf{110}(2003), 55-76.
\bibitem{Montgomery}R. Montgomery, \textit{A Tour of Subriemannian Geometry, Their Geodesics and
Applications},
Mathematical Surveys and Monographs, vol. 91, 2002.
\bibitem{Pa} S. D. Pauls, Nonparametric minimal surfaces in the
Heisenberg group, to appear in Geometriae Dedicata.
\bibitem{Pa1}S.D. Pauls, A notion of rectifiability modelled on Carnot groups,
preprint 2001.
\bibitem{Sch}J. A. Schouten, On nonholonomic connections, Proc.
Amsterdam Nederl. Akad. Wetensch., Ser. A, \textbf{31}(1928),
291-299, Jbuch 54, 758.
\bibitem{Sussmann}H. J. Sussmann, Orbits of families of vector
fields and integrability of distributions, Trans. Amer. Math. Soc.
\textbf{80}(1973), 171-188.
\bibitem{TY}K. H. Tan, X. P. Yang, Horizontal connection and horizontal mean curvature in Carnot
groups, preprint, 2003.
\bibitem{Vagner1}V. V. Vagner, Differential geometry of
nonholonomic manifolds, Tr. Semin. Vectorn. Tenzorn. Anal.
\textbf{213}(1935), 269-314(Russian).
\bibitem{Vagner2}V. V. Vagner, Differential geometry of
nonholonomic manifolds, In: \textit{The VIII-th International
Competition for the N. I. Lobatschewski Prize Report.} The Kazan
Physico-Mathematical Society, Kazan, 1940.
\bibitem{Vagner3}V. V. Vagner, A geometric interpretation of
nonholonomic dynamical systems, Tr. Semin. Vectorn. Tenzorn. Anal.
\textbf{5}(1941), 301-327.
\bibitem{Vagner4}V. V. Vagner, Geometria del calcolo delle
variazioni, vol. 2 C.I.M.E., Roma.
\bibitem{VF}A. M. Vershik, L. D. Faddeev, Lagrange mechanics in an
invariant form, In: Probl. Theor. Phys. Vol.II, 129-141. English
transl.: Sel. Math. Sov., \textbf{1}(1981), 339-350.
\bibitem{VG} A. M. Vershik, V. Ya. Gershkovich, Nonholomic
problems and the theory of distributions, Acta. Appl. Math.,
\textbf{12}, 181-209, 1988.
\end{thebibliography}
\end{document}